\documentstyle[12pt]{article}
\newcommand{\const}{\mathop{\rm const}\limits}

\newcommand{\mes}{\mathop{\rm mes}\limits}

\newcommand{\cov}{\mathop{\rm cov}\limits}

\newcommand{\vraisup}{\mathop{\rm vraisup}\limits}

\newcommand{\rank}{\mathop{\rm rank}\limits}

\newcommand{\mod}{\mathop{\rm mod}\limits}

\newcommand{\trace}{\mathop{\rm trace}\limits}

\newcommand{\sign}{\mathop{\rm sign}\limits}

\newcommand{\Law}{\mathop{\rm Law}\limits}

\begin{document}

\begin{center}

{\bf THEORY OF APPROXIMATION AND  \\
CONTINUITY OF RANDOM PROCESSES.} \par

\vspace{4mm}

 $ {\bf E.Ostrovsky^a, \ \ L.Sirota^b } $ \\

\vspace{4mm}

$ ^a $ Corresponding Author. Department of Mathematics and computer science, Bar-Ilan University, 84105, Ramat Gan, Israel.\\
\end{center}
E - mail: \ galo@list.ru \  eugostrovsky@list.ru\\
\begin{center}
$ ^b $  Department of Mathematics and computer science. Bar-Ilan University,
84105, Ramat Gan, Israel.\\

E - mail: \ sirota3@bezeqint.net\\

\vspace{3mm}
                    {\sc Abstract.}\\

 \end{center}

 \vspace{4mm}

 We  introduce and investigate a new  notion of the theory of approximation-the so-called {\it  degenerate approximation, }
i.e. approximation of the function of two (and more) variables (kernel)  by means of degenerate function (kernel).\par
 We apply obtained results to the  investigation of the local structure of random processes, for example, we find
 the necessary and sufficient condition for continuity of Gaussian and non-Gaussian  processes, some conditions for weak
 compactness and convergence of a family of random processes, in particular, for Central Limit Theorem in the space of
 continuous functions.\par
 We give also many examples in order to illustrate the exactness of proved theorems. \par
  \vspace{4mm}

\par

\vspace{4mm}

 \par

\vspace{4mm}

\section{Notations. Statement of problem.}

\vspace{3mm}

 Let $  T = [0, 2 \pi], \ T^d = [0, 2 \pi]^d; \ $  for integrable $ 2 \pi $ periodical  on all variables
 function $ f: T^d \to R $  we denote as ordinary

$$
||f||_{p,d} = ||f||_{p} = \left[\int_{T^d} |f(x)|^p \ dx   \right]^{1/p}, \ 1 \le p < \infty;
$$

$$
||f|| = ||f||_{\infty,d} = ||f||_{\infty} = \vraisup_{t \in T^d} |f(t)|;
$$

$$
d = 1 \ \Rightarrow  c_k = c_k(f) := (2 \pi)^{-1} \int_0^{2 \pi} f(x) \ e^{-i k x} \ dx;
$$

$$
\omega_p(f,\delta) := \sup_{h: ||h|| \le \delta} || f(\cdot + h) - f(\cdot)||_p, \ \omega(f,\delta) =
\omega_{\infty}(f,\delta), \ \delta \in [0,2 \pi]).
$$

 Let also $ n = \vec{n} = \{ n_1, n_2,\ldots, n_d \}, \ n_j = 0,1,2, \ldots $ be a multiindex; we denote by $ A(\vec{n}) = A(n) $
the set of all multivariate trigonometrical polynomials of a degree less or equal $ n: \ A(\vec{n}) = \{ g = g(x), x \in T^d \}, $

$$
g(x) = \sum_{\vec{k} \le \vec{n} } a(\vec{k}) \ e^{ i \vec{k} \cdot \vec{x}  },
$$
where the inequality $ \vec{k} \le \vec{n}  $ is understood coordinate-wise: $ \forall j = 1,2,\ldots,d \ k_j \le n_j. $ \par

 The following notion is very important in the approximation theory:  an error $ E(\vec{n},f)_p $ of the best
trigonometrical approximation of a function$ f(\cdot) $ in the $  L_p $ sense:

$$
E(\vec{n},f)_p \stackrel{def}{=} \inf_{ g \in A(\vec{n})} ||f - g||_p, \ 1 \le p < \infty; \
$$

$$
E(\vec{n},f) = E(\vec{n},f)_{\infty}   = \inf_{ g \in A(\vec{n})} \sup_x | f(x) - g(x)|.
$$
 See e.g. the classical monographs \cite{Achieser1}, \cite{DeVore1}, \cite{Korneichuk1}, \cite{Lorentz1}, \cite{Natanson1},
\cite{Rivlin1}, \cite{Timan1} etc.\par
  For instance, there are obtained Jackson's type estimations of a view

$$
E(n,f)_p \le C(r,p) \ n^{-r} \ \omega(f^{(r)}, \pi/n)_{p}, \ r=0,1,2, \ldots, \ p \in [1,\infty], \ d = 1
$$
with exact value of "constants" $ C(r,p), $  see \cite{Korneichuk1}, \cite{Natanson1}, \cite{Timan1}.\par
 Here  and further  $  \omega(f, \delta)_p, \ \delta \in (0, 2 \pi), \ 1 \le p \le \infty $ will denote  the
{\it one-dimensional} ordinary $ L_p $ modulus of continuity  of the function $ f; $ \par

$$
\omega(f, \delta) :=  \omega(f, \delta)_{\infty} \stackrel{def}{=} \sup_{0 < |h| < \delta} \sup_t |f(t + h) - f(t)|.
$$

\vspace{3mm}

{\bf We  introduce and investigate in this article the  degenerate approximation,
 and apply the obtained results to the theory of random fields, for instance for  finding  of
 necessary and sufficient condition for continuity of Gaussian and non-Gaussian  processes. } \par

\vspace{3mm}

 The paper is organized as follows. In the next section we define and investigate a new notion: degenerate
approximation of a function of several (two, for beginning) variables by means of linear combinations products
of a one variable functions.  The third section contains the main results: necessary and sufficient condition
for continuity of Gaussian  random processes  in the terms of   degenerate approximation of its covariation function. \par
  We obtain in the 4th section on the basis of main result some new sufficient conditions for continuity of
Gaussian processes, not necessary to be stationary. The next section  is devoted to generalization for
non-Gaussian case, the 6th section we ensue some sufficient conditions for continuity of non-Gaussian processes.\par
 In the seventh section we formulate and prove the necessary and sufficient condition for convergence of
the sequence of random variables with probability one.\par
 In the last section there are several concluding remarks. \par

\vspace{3mm}
\section{ Degenerate approximation. }

\vspace{3mm}

{\bf Definition 2.1.}  The function  of two variables $ R = R(t,s), \ t,s \in T $  of a view

$$
R(t,s) = \sum_{k=1}^{n_1} \sum_{l=1}^{n_2} a_{k,l} \phi_k(t) \psi_l(s), \ a_{k,l} = \const \eqno(2.1)
$$
is said to be {\it degenerate;}  more exactly, linear degenerate. \par
 {\it  It will be presumed that if all the functions} $ h_1(t) = R(t,s) $ {\it and } $  h_2(s) = R(t,s) $
 {\it  belong to some Banach spaces, then all the functions  } $  \{ \phi_k(\cdot)  \} $
 {\it and } $ \psi_l(\cdot)  $ {\it belong to at the same space.  }\par
  For instance, if

  $$
  \forall s \in T \ \Rightarrow \int_T  |R(t,s)|^p \ dt < \infty,
  $$
 then $ \phi_k(\cdot) \in L_p. $\par
  If the function $ (t,s) \to R(t,s) $ is  continuous, then all the functions $ \phi_k(t), \ \psi_l(s)  $ are continuous. \par

 We can and will assume that both the systems of a functions $ \{\phi_k(\cdot) \}, \ k = 1,2,\ldots, n_1 $ and
$ \psi_l(\cdot), \ l=1,2, \ldots, n_2  $ are linear independent and moreover orthonormal:

$$
\int_T \phi_k(t) \phi_l(t) \ dt = \delta_{k,l} = \int_T \psi_k(t) \psi_l(t) \ dt,
$$
where $ \delta_{k,l}  $ is Kroneker's symbol. \par
 We will call the equality (2.1) as {\it representation} of a function $ R(t,s) $ and the pair of numbers
 $ (n_1,n_2) $ will named $ \rank(R): \ (n_1,n_2) := \rank(R). $ \par

  \vspace{4mm}

{\bf Definition 2.2.} Suppose the degenerate function $ R(t,s) $ is symmetrical: $  R(t,s) = R(s,t). $
The representation of a view

$$
R(t,s) = \sum_{k=1}^{n} b_{k} \phi_k(t) \phi_k(s), \ b_{k} = \const \eqno(2.2)
$$
is called {\it symmetrical}.\par

\vspace{4mm}

 The set all of degenerate functions  with finite rank $ (n_1, n_2) $ will be denotes by
$ Q(n_1,n_2); $ the set all of symmetrical degenerate functions  with finite rank $ (n, n) $ will be denotes by
$ Q^{(s)}(n). $ \par

\vspace{4mm}

 This notion is very important and may be used, for instance, in the theory of integral equations, where
in the case when the kernel of equation $ R(t,s):  $

$$
X(t) = \int_T R(t,s) \ X(s) \ ds = f(t)
$$
is degenerate, then this equation may be completely investigated and solved in explicit view.\par

\vspace{3mm}

 Analogously, in the game theory, where $ R(t, s)  $ is the so-called pay function in the game between two
 players  with zero sum, when $  R(\cdot, \cdot)  $ is degenerate, then this game may be solved
in explicit view.

\vspace{3mm}

 The application of this notion in the theory of random processes,  where $ R(t,s) $ is covariation function
for  Gaussian random process (field), obviously symmetrical and non-negative defined,  will be discussed in the
third section. \par

\vspace{3mm}

In the more general statement belonging to A.N.Kolmogorov (see \cite{Vitushkin1} )
considered in this article problem  is called
"Separability of multivariate functions"; see e.g. some recent   publications
\cite{Goda1}, \cite{Kuo1}, \cite{Sobol1}; where are considered in particular applications in the
numerical methods and in statistics.\par

\vspace{4mm}

{\bf Definition 2.3. \ Error of degenerate approximation. } \par
   Let $ R: T^2 \to R^1 $ be integrable function and let $ n_1, n_2  $ be natural numbers. We define the value
 $ D(n_1, n_2;R)_p $ as an error of minimal $ L_p $ approximation of the function $  R(\cdot, \cdot) $ by means of
 degenerate functions of the rank $ (n_1, n_2): $

 $$
 D(n_1, n_2;R)_p  = \inf_{ G \in Q(n_1,n_2)} || R - G||L_p(T^2). \eqno(2.3)
 $$
 \vspace{4mm}
 As before, we define for the continuous function $ R(\cdot, \cdot) $

 $$
 D(n_1, n_2;R)  = D(n_1, n_2;R)_{\infty}  = \inf_{ G \in Q(n_1,n_2)} \sup_{t,s} | R(t,s) - G(t,s)|. \eqno(2.4)
 $$
\vspace{4mm}
  Analogously we define for the symmetrical functions $ R(t,s) $

 $$
 D^{(s)}(n;R)_p  = \inf_{ G \in Q^{(s)}(n)} || R - G||L_p(T^2). \eqno(2.5)
 $$

\vspace{4mm}

{\bf Example 2.1.}  Let $  R = R(t,s) $ be non-zero continuous non-negative  definite function.
for example, the covariation function of Gaussian mean-square continuous process. Then
the value $ D^{(s)}(n;R)_2 $ and the optimal approximate degenerate function $ Q(n; R) $ may be calculated
in explicit view. Namely, consider the following equation for eigen  values and functions:

$$
\int_T R(t,s) \ \phi_k(s) \ ds = \lambda_k \ \phi_k(t). \eqno(2.6)
$$
 It follows from the theorem of Hilbert-Schmidt that the equation (2.6) has an enumerable set of eigen
values $ \lambda_k, \ k=1,2,\ldots $ and  orthonormal complete in $ L_2(T) $ sense sequence of eigen (continuous)
functions    $  \{ \phi_k \}. $\par
  It follows from the Mercer's theorem that

$$
R(t,s) = \sum_k \lambda_k \ \phi_k(t) \ \phi_k(s), \eqno(2.7)
$$
and the series (2.7) converges uniformly. \par

 Moreover,

 $$
 \sum_{k=1}^{\infty} \lambda_k = \trace(R)  = \int_T R(t,t) \ dt < \infty.
 $$

 Since the kernel $ R(t,s) $ is non-negative definite, all the values $ \{\lambda_k\} $ are non-negative; we can restrict
ourselves only  by positive eigen values $ \lambda_k; \lambda_k > 0. $\par

 We can suppose without loss of generality that the sequence of the eigen values decreases:

$$
\lambda_1 \ge \lambda_2 \ge \lambda_3 \ge \ldots \ge  \lambda_n \ge \ldots; \eqno(2.8)
$$
the case when only finite set of eigen values is positive is trivial: it follows from (2.7) that $ R(t,s) $ is
symmetrical degenerate.\par

 It is easy to calculate that the optimal approximate degenerate function $ Q(n; R) $ for the kernel $ R(\cdot, \cdot) $
is follows:

$$
 Q(n; R)(t,s) = \sum_{k=1}^n \lambda_k \ \phi_k(t) \ \phi_k(s), \eqno(2.9)
$$
and correspondingly

$$
D^{(s)}(n;R)_1 = \sum_{k=n+1}^{\infty} \lambda_k  \to 0, \ n \to \infty. \eqno(2.10)
$$

\vspace{4mm}

{\bf Example 2.2.} Let in the example (2.1)  for some integer value $ n \ \lambda_n = \lambda_{n+1}.  $  Then in the capacity
of  the optimal approximate by degenerate function of the finite fixed rank $ n $
$ Q(n; R) $ may be choose might as well the other function

$$
 Q(n; R)(t,s) = \sum_{k=1}^{n-1} \lambda_k \ \phi_k(t) \ \phi_k(s)  + \lambda_n \phi_{n+1}(t) \ \phi_{n+1}(s). \eqno(2.11)
$$
 This reason imply  that the optimal approximate degenerate function may be not unique. \par

\vspace{4mm}
 The main result of this section is follows.\par
 \vspace{4mm}

 {\bf Theorem 2.1.}

$$
 D(2 n_1, \ 2 n_2; R)_p \le 1 \cdot E(n_1,n_2; R)_p, \ 1 \le p \le \infty,  \eqno(2.12)
$$
{\it  where for all the values $ p $  the constant   "1" in (2.12) is the best possible}.\par
\vspace{4mm}

{\bf Proof. \ A. Inequality.}\par
 The upper estimate

$$
 D(2n_1, \ 2n_2; R)_p \le E(n_1,n_2; R)_p, \ 1 \le p \le \infty,
$$
is very simple. Indeed, there exists a two-variate  trigonometrical polynomial $ Q(n_1, n_2)(t,s) $ of degree $  n_1,n_2,  $
not necessary to be unique, for which

$$
|| \ R -   Q(n_1, n_2)   \ ||_p = E(n_1,n_2; R)_p.
$$
 It remains to note that the function $ Q(n_1, n_2) (t,s)  $ is degenerate with the rank $ (2n_1, \ 2n_2). $\par

\vspace{4mm}

{\bf Proof. \ B. Exactness.} \par
  It is sufficient to consider the case $ n_1 = n_2 = n $ and symmetrical function $ R(\cdot, \cdot).  $ \par
 We need firs of all to specify the statement of this problem. Let us define

 $$
 V = V(p) = \overline{\lim}_{n \to \infty} \sup_{R \in L_p(T^2)} \left[ \frac{D^{(s)}(2n;R)_p} {E(n,n,R)_p} \right]. \eqno(2.13)
 $$
 It follows from the upper estimate (2.12) that $  V \le 1; $ it remains to prove opposite inequality.\par
  We consider as an example the following lacunar series:

$$
R(t,s) = R_l(t,s) = \sum_{k=1} ^{\infty} a_k \ \cos(n_k t) \ \cos(n_k s), \eqno(2.14)
$$
where $ a_k > 0, \ \sum_k a_k < \infty, \ n_k  $ be a increasing sequence of integer numbers such that

$$
\frac{n_{k+1}}{n_k} = 2 p_k  + 1, \
$$
$  p_k  $ is arbitrary sequence of integer numbers under condition $ p_k \ge 2. $ \par
 Obviously, the function $ R_l(t,s) $ is continuous since the series in (2.14) converges uniformly.\par
 It is known from the theory of lacunar trigonometrical  series, see, for example, \cite{Timan1}, chapter 8,
that

$$
E(n,n,R_l)_p  = \sum_{k = \nu + 1}^{\infty} a_k, \ n_{\nu} \le n < n_{\nu + 1}. \eqno(2.15)
$$
 On the other hands, it follows from the example 2.1 that

 $$
 D(2n_{\nu}, \ 2n_{\nu}; R_l)_p =  \sum_{k = \nu + 1}^{\infty} a_k, \ n_{\nu} \le n < n_{\nu + 1}.
 $$
 Therefore $  V(p) \ge 1, $ Q.E.D. \par
\vspace{4mm}

{\bf Example 2.3.} Let $ R(t,s) $ be degenerate function of finite rank but not a trigonometrical polynomial,
for instance

$$
R(t,s) = \phi(t) \phi(s), \ \phi \notin \cup_{n=1}^{\infty} A(n).
$$
Then $ Q(n,n; R) = 0 $ for all sufficiently greatest values $ n, $ but for all the values $ n \  E(n,n; R) > 0. $

\vspace{4mm}
\section{ Main results: criterion for continuity of \\
 Gaussian  random processes  in the terms  \\
  of  degenerate approximation. }

\vspace{4mm}

 \begin{center}

{\bf A. Statement of problem.} \par

\end{center}

\vspace{3mm}

Let $ \xi(t), \ t \in T  $ be $ 2\pi $ periodical separable random process (r.p.)
defined aside from the set $ T $ on some probabilistic space with probabilistic measure $ {\bf P }, $
expectation $ {\bf E}  $ and covariation $ \cov. $ Question: under what conditions the  r.p. $ \xi(t) $
is continuous with probability one:

$$
  {\bf P} (\xi(\cdot) \in C(T)) = 1? \eqno(3.0)
$$

 This problem has a long history. The firs result belongs to A.N.Kolmogorov  \cite{Kolmogorov1} and
E.E.Slutsky \cite{Slutsky1}.  For the separable random fields $ \xi(\vec{t}), \ \vec{t} \in [0,1]^d  $ the correspondent
result belongs to  V.Vinkler  \cite{Vinkler1}: if

$$
{\bf E} |\xi(\vec{t}) -  \xi(\vec{s})|^{\alpha} \le C \cdot ||\vec{t} - \vec{s} ||^{d + \beta},  C < \infty,
0 < \alpha, \beta < \infty,
$$
then $ {\bf P} (\xi(\cdot) \in C([0,1]^d) = 1. $\par
  Further results was obtained by means of the so-called entropy terms introduced by
 R.M.Dudley \cite{Dudley1} and X.Fernique \cite{Fernique1} -  \cite{Fernique3}. \par
  A very interest results for Gaussian {\it stationary} processes and fields belong to N.Nisio \cite{Nisio1} and
 H.Watanabe  \cite{Watanabe1}. \par
   Recently appear many works based on the so-called notion of {\it majorizing and minorizing measures,} see e.g. the articles and books of
X.Fernique \cite{Fernique1}, M.Ledoux and Talagrand \cite{Ledoux1},  M.B.Marcus and L.A.Shepp  \cite{Marcus1},
M.Talagrand  \cite{Talagrand1} - \cite{Talagrand5},  W.Bednorz
\cite{Bednorz1}-  \cite{Bednorz4};  see also  \cite{Garsia1}, \cite{Heinkel1},  \cite{Kwapien1},  \cite{Ostrovsky1},
\cite{Ostrovsky100},  \cite{Ostrovsky101},  \cite{Ostrovsky102},  \cite{Ostrovsky103},  \cite{Ostrovsky207} - \cite{Ostrovsky210},
\cite{Ral'chenko1}. \par
 {\it We intend in this and in the next sections to find some conditions (necessary and sufficient,  sufficient) for the equality (3.0) in the
terms of approximation theory, in particular, in the terms of degenerate approximation.} \par

\vspace{4mm}

\begin{center}

{\bf B. Necessary and sufficient condition for continuity of \\
  Gaussian process.}\par

\end{center}

 Note that this problem was solved in the another terms: by means of Jackson's polynomials in \cite{Ostrovsky208},
in the geometrical notions of Hilbert's space  generated by the Gaussian r.p. in \cite{Sudakov1},  in the terms of the so-called
partition schemes (generic  chaining) - in \cite{Talagrand2}, chapters 3,4.  For the  stationary Gaussian processes the famous criterion for
continuity  was obtained in the entropy terms by X.Fernique in \cite{Fernique1}. \par

\vspace{4mm}

 We consider in this subsection the centered Gaussian periodical process $ \xi(t), \ t \in T, \ {\bf E} \xi(t) = 0, $
with (non-negative definite)  covariation function  $ R(t,s) = \cov(\xi(t), \xi(s)) = {\bf E} \xi(t) \xi(s).  $ \par

 Let $ \{ \eta_k \}, \ k=1,2,\ldots $
be a sequence of independent standard distributed Gaussian r.v.  Let also as  in (2.7), (2.8) $ \lambda_k, \ \phi_k(t) $
be  correspondingly a sequence of eigen values and normed eigen function for the kernel  $ R. $ \par
 We define the following family of a semi-norms:

$$
m \ge n+1 \ \Rightarrow \tau_n^m(R) \stackrel{def}{=} {\bf E} || \sum_{k= n+1}^m \sqrt{\lambda_k} \ \eta_k \ \phi_k||_{\infty} =
$$

$$
(2 \pi)^{(n-m)/2 } \int_{R^{m-n}}  \exp \left( -0.5 \sum_{k=n+1}^m x_k^2 \  \right)\cdot
|| \sum_{k=n+1}^{m} x_k \sqrt{\lambda_k} \ \phi_k(\cdot) ||_{\infty}  \cdot \prod_{k=n+1}^m dx_k; \eqno(3.1)
$$

$$
\tau(R) := \overline{\lim}_{n \to \infty} \sup_{m \ge n+1} \tau_n^m(R),
$$
where in the case when $ \tau_n^m(R) $ is not defined, we  take by definition $\tau_n^m(R)  = + \infty. $ \par

\vspace{3mm}

{\bf Theorem 3.1.}  {\it  In order to the separable $ 2 \pi $ periodical Gaussian r.p. } $ \xi(t) $ {\it   has continuous with
probability one trajectories, is necessary and sufficient that its covariation function  } $ R(t,s) $ {\it is  continuous and}

$$
\overline{\lim}_{n \to \infty} \sup_{m \ge n+1} \tau_n^m(R) = 0, \eqno(3.2)
$$
{\it or equally} $  R \in \ker(\tau). $ \par

\vspace{3mm}

{\bf Proof. Sufficiency.} Let the  conditions of theorem 3.1.be satisfied; then we can apply Mercer's theorem (2.7) with
the order (2.8). We can write the Karunen-Loev expression for $ \xi(\cdot): $

$$
\xi(t) = \sum_{k=1}^{\infty} \sqrt{\lambda_k} \ \zeta_k \ \phi_k(t). \eqno(3.3)
$$
 Here $  \{ \zeta_k \} $ is a sequence of independent standard distributed Gaussian r.v.
 See in detail  \cite{Ostrovsky209}; in particular, it is prove therein that the eigen function $ \phi_k(t) $  are continuous.  \par
 Denote the partial sum of Karunen-Loev expression as $  \xi_n = \xi_n(t): $

$$
\xi_n(t) = \sum_{k=1}^n  \sqrt{\lambda_k} \ \zeta_k \ \phi_k(t).
$$
 Let $  n(l), \ l=1,2,\ldots $ be any strictly increasing  subsequence of natural numbers such that

$$
\tau_{n(l)+1}^{n(l+1) } \le 2^{-l},  \ l \ge l_0 = \const,
$$
and define

$$
S_l = S_l(t) = \sum_{k=n(l) + 1}^{n(l+1)} \sqrt{\lambda_k} \ \zeta_k \  \phi_k(t) = \xi_{n(l+1)}(t) - \xi_{n(l)}(t).
$$
then

$$
{\bf E} ||S_l||_{\infty} = \tau_{n(l)+1}^{n(l+1) } \le 2^{-l}, \
\sum_{l \ge l_0}  {\bf E} ||S_l||_{\infty} = \sum_{l \ge l_0} \tau_{n(l)+1}^{n(l+1) }  < \infty,
$$
and following the series

$$
\sum_{l \ge 1} ||S_l||_{\infty}
$$
converges with probability one.
Therefore, its partial sums, i.e. the subsequence of r.p. $ \xi_{n(l)}(t) $  converges uniformly with probability one.
Thus, the r.p. $ \xi(t) $ is continuous $ ( \mod{\bf P}) $. \par

\vspace{3mm}
{\bf Proof. Necessity.} It can be assumed here that the centered Gaussian r.p. $ \xi(t) $ is continuous a.e. Therefore,
 its covariation function $ R(t,s) $ is continuous and consequently the Karunen-Loev expression there exists. \par
 It follows from theorem of K.Ito and M. Nisio \cite{Ito1} that in considered case  the Karunen-Loev expression
converges uniformly with probability one:

$$
{\bf P} ( \lim_{n \to \infty} || \ \xi(\cdot) - \xi_n(\cdot) \ ||_{\infty} = 0) = 1.  \eqno(3.4)
$$
 In turn, we conclude based on the equality (3.4) that

 $$
 \lim_{n \to \infty} {\bf E} || \ \xi(\cdot) - \xi_n(\cdot) \ ||_{\infty} = 0,   \eqno(3.5)
$$
see  \cite{Ostrovsky211}. \par
 We get ultimately using triangle inequality

 $$
 \tau_n^m(R) = {\bf E} || \xi_m(\cdot) - \xi_n(\cdot) ||_{\infty} \le {\bf E} || \xi_m(\cdot) - \xi(\cdot) ||_{\infty} +
 $$

$$
{\bf E} || \xi(\cdot) - \xi_n(\cdot) ||_{\infty} \to 0, \ n,m \to \infty.
$$
 This completes the proof of theorem 3.1. \par

\vspace{3mm}

 \begin{center}

 {\bf C. Some comments.} \\

 \end{center}

\vspace{3mm}

{\bf 1.  Non-centered case.}  Let $  \xi(t) $ be separable Gaussian distributed r.p. with non-zero expectation
 $ {\bf E} \xi(t) = a(t),  \ t \in T. $  The r.p. $ \xi(t) $ is continuous iff $ a(t) $ is continuous and the
 centered Gaussian process $ \xi^{(0)}(t) := \xi(t) - a(t) $ satisfies the condition of theorem 3.1. \par

\vspace{3mm}

{\bf 2. Discontinuous case.}  If the condition of theorem 1.3 are non satisfied, then

$$
{\bf P} (\xi(\cdot) \in C(T) ) = 0,
$$
(zero-one law); see \cite{Fernique1}, \cite{Ostrovsky209}. \par

\vspace{3mm}

{\bf 3. Non-periodical case.}  The condition of periodicity  $ \xi(t \pm 2 \pi) = \xi(t) $ is not essential
restriction. Let for instance $ \xi(t)  $ be centered Gaussian process, $ t \in [0,1]. $  Consider the linear
extrapolation (spline)  $ \tilde{\xi}(t) $ of $ \xi(t) $ on the set $  [1, 2 \pi], $ i.e. such that

$$
\tilde{\xi}(t) = \xi(t), \ t \in[0, 1]; \ \tilde{\xi}(2 \pi) = \xi(0)
$$
and $ \tilde{\xi}(t) $  is periodical continuation  on the whole axis. \par
 The Gaussian r.p. $ \xi(\cdot) $ is continuous a.e. iff the periodical Gaussian r.p. $ \tilde{\xi}(\cdot) $
is continuous. \par

\vspace{3mm}

{\bf 4. Other norms.} The expression for the semi-norm $  \tau_n^m (R) $ may be rewritten as follows:

$$
\tau_n^m (R) =  || \ \max_{t \in T} | \sum_{k= n+1}^m \sqrt{\lambda_k} \ \eta_k \ \phi_k(t) | \  ||L_1(\Omega, {\bf P}). \eqno(3.6)
$$
 But  instead the $ L_1(\Omega) $ norm may be used more strong rearrangement invariant
 norm, up to Orlicz's  norm over our probabilistic space with $ N - $  function

$$
N(u) = \exp(u^2/2) - 1;
$$
see \cite{Ostrovsky211}.\par

\vspace{3mm}

{\bf 5. Banach space valued Gaussian r.v.}  Let  $ B $ be separable Banach space equipped with the norm
$ || \cdot ||B $ and $  \xi $  be centered weak  Gaussian
distributed r.v. with covariation operator $ R. $  We ask: under some condition

$$
{\bf P} (\xi \in B) = 1?  \eqno(3.7)
$$

 It is known \cite{Kwapien1}, \cite{Wakhaniya1}, chapter 5, section 5  that if (3.7) there holds, then
there exists a sequence  $  \{ g_n \}  $ of topological free non-random elements of  the space $  B  $  and a sequence  of
independent standard distributed Gaussian r.v. $  \{  \zeta_n \} $ such that

$$
\sum_n ||g_n||^2 B  < \infty
$$
and

$$
\xi = \sum_{n=1}^{\infty} \zeta_n \ \ g_n. \eqno(3.8)
$$
 Define  for natural numbers  $ n,m: \  m \ge n+1  $ the following semi-norm on the space of all symmetrical
operators $  \{   R \} $

$$
\tau_n^m(R)[B] = {\bf E} || \sum_{k=n+1}^{m} \zeta_k \ \ g_k  ||B =
$$

$$
(2 \pi)^{(n-m)/2 } \int_{R^{m-n}}  \exp \left( -0.5 \sum_{k=n+1}^m x_k^2 \  \right) \cdot
|| \ \sum_{k=n+1}^{m} x_k g_k \ ||B   \cdot \prod_{k=n+1}^m dx_k. \eqno(3.9)
$$
 We conclude as before that the series (3.8) converges in the norm $ || \cdot||B $ with probability one
and  as a consequence  (3.7) there holds if and only if

$$
\overline{\lim}_{n \to \infty} \sup_{m \ge n+1} \tau_n^m(R)[B] = 0. \eqno(3.10)
$$

  Obviously, instead the $ L_1(\Omega) $ norm may be used more strong norm, up to Orlicz's norm over our probabilistic space
with $ N - $  function $ N(u) = \exp(u^2/2) - 1. $ \par
 Taking for instance the Hilbertian norm induced by $ N - $ function of a view
$  N(u) = u^2, $ we get the classical criterion for equality (3.7) for the case of separable Hilbert space $  H: $

$$
{\bf P} (\xi \in H) = 1 \  \Leftrightarrow \ \trace(R) < \infty.
$$

\vspace{3mm}

{\bf 6.  Estimation of tail of distribution of } $ ||\xi||B. $  We have under conditions of last subsection the
following estimate:

$$
 E ||\xi||B \le \sup_n \sup_{m \ge n + 1} \tau_n^m(R)[B] =:  \overline{\tau} < \infty. \eqno(3.11)
$$
 Further, in the classical book of N.N.Wakhaniya, W.I.Tarieladze and S.A.Chobanjan  \cite{Wakhaniya1}
on the pages 263-264 is proved  the following important inequality for the arbitrary centered Gaussian r.v.
$ \xi $ in the space $ B: $

$$
 \left[ {\bf E} (||\xi||B)^p \right]^{1/p}  \le c(p,s) \ \left[ {\bf E} (||\xi||B)^s \right]^{1/s}.
$$
 In particular, there exists an absolute finite positive constant $ C_1 $ such that $ c(p,1)
 \le C_1 \cdot \sqrt{p}, \ p \ge 2. $ Therefore

$$
 \left[ {\bf E} (||\xi||B)^p \right]^{1/p} \le C_1 \cdot \overline{\tau} \cdot \sqrt{p}
$$
or equally

$$
{\bf P} (||\xi||B > u)  \le \exp \left( - C_2 (u/ \overline{\tau})^2  \right),  \  u > \overline{\tau}. \eqno(3.12)
$$

\vspace{3mm}

{\bf 7. Special Banach spaces.}  The case when the Banach space $  B  $ does not contains the
subspace isomorphic to the space $ c_0 $ is considered in the monograph of V.V.Buldygin
\cite{Buldygin1}, p. 116-127; see also reference therein. \par

\vspace{3mm}

\begin{center}

{\bf D. Boundedness of the Gaussian random processes. } \\

\end{center}

 The Gaussian (or with other distribution) r.p. $ \xi(t) $ is called  bounded (almost surely) if

 $$
 {\bf P} (\sup_t |\xi(t)| < \infty )  = 1. \eqno(3.13)
 $$

 The criterion for the boundedness of the separable centered Gaussian periodical r.p.  $  \xi(t) $ with continuous
covariation  function $  R(\cdot, \cdot) $  is follows:

$$
\tilde{\tau}(R) \stackrel{def}{=} \sup_n \sup_{m \ge n+1} \tau_n^m (R) < \infty; \eqno(3.14)
$$
we retain the previous notations. \par
 A.N.Kolmogorov formulated  as a hypotheses  and Yu.K. Belyaev proved in \cite{Belyaev1}  the following
alternative  for {\it  stationary }  Gaussian centered separable r.p. $ \xi(t): $ either it is continuous
a.e. or is unbounded on each non-empty interval. \par
 Note that this alternative is not true for the non-stationary Gaussian processes.  For example: let
 $  \zeta(t) = \epsilon \cdot \sign(t), \ \Law(\epsilon) = N(0,1),  \ t \in [-1,1];  $
 then the Gaussian r.p. $ \zeta(t)  $ is bounded and discontinuous. \par
   But this process has discontinuous covariation function. We offer another example of a centered discontinuous {\it bounded}
Gaussian process $ \beta(t)  $ with continuous  covariation function. Let $ w(t), \ t \in [0,e^{-4}]  $ be ordinary
Brownian motion; we define $ \beta(t)  $ as follows:

$$
 \beta(t) = \frac{w(t)}{ \sqrt{2\ t \ \log | \log t|}}, \ t > 0,  \eqno(3.15)
$$
and $ \beta(0) = 0. $   Evidently, $ \beta(t) $ is centered, Gaussian  and mean square continuous. \par
 It follows from the Law of Iterated Logarithm  (LIL) that

 $$
 \overline{\lim}_{t \to 0+}  \beta(t) = 1 = - \underline{\lim}_{t \to 0+}  \beta(t),  \eqno(3.16)
 $$
therefore $ \beta(t) $ is bounded and discontinuous.\par

\vspace{4mm}

\section{ Some new sufficient conditions for continuity of Gaussian processes.}

\vspace{4mm}

 The author be aware that the conditions of theorem 3.1 are hard to verify. Further in this section
we will obtain some simple {\it sufficient}  conditions for continuity of (periodical) separable Gaussian
centered r.p. $ \xi(t)  $ based on some estimations for important for us semi-norm $ \tau_n^m(R). $\par

\vspace{4mm}

{\bf Theorem 4.1.}  {\it Assume that there exists a strictly increasing non-random  sequence of natural numbers }
$  \{  n(k) \}, \ k=1,2,\ldots $ {\it such that} $ n(1) = 1, $

$$
\Sigma_E \stackrel{def}{=}  \sum_{k=1}^{\infty}  E^{1/2}([n(k)/2], [n(k)/2]; R)  \cdot \sqrt{\log n(k+1)} < \infty. \eqno(4.0)
$$
{\it Then} $ \ {\bf P} (\xi(\cdot) \in C(T) ) = 1 $  {\it and moreover for some absolute constants }  $ C_3, C_4 $

$$
 \left[ {\bf E} (||\xi||B)^p \right]^{1/p} \le C_3 \cdot \Sigma_E \cdot \sqrt{p}, \ p \ge 1,
$$

$$
{\bf P} (||\xi||B > u)  \le \exp \left( - C_4 (u/\Sigma_E)^2  \right),  \  u > \Sigma_E.
$$

\vspace{3mm}

{\bf Proof.} \par

\vspace{3mm}

{\bf 1.} We intend   to apply in order to obtain the good degenerate approximation for the kernel $ R(\cdot,\cdot) $
the well-known  Vallee-Poussin sums \cite{Poussin1}, \cite{Timan1}, chapter 5.
 Recall that the Vallee-Poussin kernel $ K_{n,p}(t) $ is defined as follows:

$$
K_{n,p}(t) = \frac{\sin((2n+1-p)t/2) \cdot \sin((p+1)t/2)}{2(p+1) \sin^2t/2}.
$$
 It is known that $ K_{n,p}(t)  $ is trigonometrical polynomial of degree  $ n. $  \par

 The Vallee-Poussin  approximation (sum) $ V_{n,p}[f](t) $ for a periodical integrable function $ f $
may be defined as follows:

$$
V_{n,p}[f](t) := [f * K_{n,p}](t)
$$
(periodical convolution). We pick hereafter for definiteness for the values $ n \ge 4  \ p = p(n) := [n/2],  $ (integer part),
so that  $ V_{n,p}[f](t)  \in A(n) $ and

$$
|| f(\cdot) - V_{n,p(n)}[f](\cdot) ||_{\infty} \le C \cdot E([n/2],f),
$$
see \cite{Natanson1}, chapter 6.  \par

\vspace{3mm}

{\bf 2. Lemma 4.1. (\cite{Ostrovsky209}).  } Let  $ \eta_n(t) $ be centered Gaussian process and simultaneously trigonometrical
polynomial of degree $  n: \ \eta_n(\cdot) \in A(n). $  Denote

$$
\sigma^2(n) = \max_t \ {\bf Var}\{ \eta_n(t) \}.
$$
 There holds:

$$
{\bf E} || \ \eta_n(\cdot) \ ||_{\infty} \le C \cdot \sigma(n) \cdot \sqrt{\log n }, \ n \ge 4. \eqno(4.1)
$$

\vspace{3mm}
{\bf 3.} Let  $ n(k), k = 1,2, \ldots  $ be arbitrary  strictly increasing sequence of natural numbers.
 Introduce the following sequence of a functions:

 $$
 Y_k (t)= [K_{n(k+1), p(n(k+1))}* \xi](t) -  [K_{n(k), p(n(k))}* \xi](t)=:
 $$
 $$
[W(n(k+1), n(k),\cdot)*\xi(\cdot)](t), \eqno(4.2)
 $$
where

$$
W(n(k+1), n(k),t)= K_{n(k+1), p(n(k+1))}(t) - K_{n(k), p(n(k))}(t).
$$

  We state: $ Y_k(t) $ is centered Gaussian process and trigonometrical polynomials with degree less than
  $ n(k+1): \  Y_k (t) \in A(n(k+1)). $ \par
   We get also applying the inequality of De la  Vallee-Poussin:

$$
\sup_t {\bf Var} [ Y_k (t)] \le C \cdot  E([n(k)/2], [n(k)/2]; R).
$$
 Applying of Lemma (4.1) yields:

 $$
 {\bf E} || \ Y_k (t) \ ||_{\infty} \le C \cdot  E^{1/2}([n(k)/2], [n(k)/2]; R)  \cdot \sqrt{\log n(k+1)}.
 $$
  We conclude on the basis of condition (4.0) that

  $$
\sum_{k=1}^{\infty}   {\bf E} || \ Y_k (t) \ ||_{\infty} < \infty,
  $$
 therefore

 $$
 \sum_{k=1}^{\infty}  || \ Y_k (t) \ ||_{\infty} < \infty \ (\mod {\bf P}). \eqno(4.3)
 $$
 Following, the subsequence

 $$
 \xi_{n(k+1)}(t) = [ \xi * V_{n(k+1)/2, p(n(k+1)/2)} ](t) = \sum_{m=1}^ k  Y_m (t)
 $$
converges uniformly in $ t; t \in T $ with probability one and hence the limit as
 $ k \to \infty $ of the r.p. $ \xi_{n(k+1)}(t) $ is continuous with probability one. \par

\vspace{4mm}

 \begin{center}

 {\bf Examples and conclusions.}

 \end{center}

\vspace{3mm}

{\bf 1.}  The condition (4.0) is satisfied, for instance, if

$$
\sum_{k=1}^{\infty} E^{1/2} \left(2^{k}, 2^{k}; R \right)  \cdot \sqrt{k}  < \infty
$$

or if

$$
\sum_{k=1}^{\infty} E^{1/2} \left(2^{2^k}, 2^{2^k}; R \right)  \cdot 2^{k/2} < \infty \eqno(4.4)
$$
 The last condition (4.4) is equivalent the following inequality:

 $$
 \sum_{k=1}^{\infty} E^{1/2}(2^{k^2}, 2^{k^2};R) < \infty.
 $$
 In turn, the condition (4.5) is satisfied, on the basis of Jackson inequality, if

 $$
 \int_0^{\infty} \omega^{1/2} \left(R, e^{-x^2} \right) dx < \infty. \eqno(4.5)
 $$
  We obtain the famous  Fernique's condition  \cite{Fernique1}. \par

\vspace{3mm}

{\bf 2.}  Let us show that  the condition (4.0)  is not necessary even for continuity of
a {\it stationary} Gaussian r.p. \par
 Consider the following example (lacunar random Fourier series):

 $$
 \xi(t) = \sum_{k=1}^{\infty} b_k \ \epsilon_k  \ \cos(2 n_k t), \eqno(4.6)
 $$
 where $  \{  \epsilon_k  \} $ is a sequence of independent standard normal distributed r.v.,
$ \ n_k  $ be a increasing sequence of integer numbers such that

$$
\frac{n_{k+1}}{n_k} = 2 p_k  + 1, \
$$
$  p_k  $ is arbitrary sequence of integer numbers under conditions $ p_k \ge 2, $
such that

$$
\lim_{k \to \infty} \frac{\log n(k+1)}{\log n(k)} = \infty.
$$

 The r.p. $  \xi(t) $ is the real part of complex centered stationary  Gaussian random process and has
a covariation function

$$
R_{\xi}(t,s) = \sum_{k=1} ^{\infty} b_k^2 \ \cos(n_k t) \ \cos(n_k s).
$$
 The necessary and sufficient condition for continuity of $ \xi(\cdot) $ is follows:

 $$
 \sum_k |b_k| < \infty. \eqno(4.7)
 $$
The  condition (4.0) of theorem 4.1 may be formulated as convergence of a series:

 $$
 \sum_k \sqrt{ \sum_{m=k}^{\infty} b^2_k} \cdot \sqrt{\log n(k+1)} < \infty. \eqno(4.8)
 $$
 Obviously, the condition (4.8) is essentially more stronger than (4.7). \par

\vspace{3mm}

{\bf 3. } Let us show that if the condition (4.5) is not satisfied, the Gaussian r.p. $ \xi(t) $
may be continuous as well as may be discontinuous.  \par
 Note firs of all that the case of non-stationary process is very simple; it is sufficient to
consider the degenerate process and covariation  $ \xi(t) = \epsilon \ \phi(t), \ R(t,s) = \phi(t) \phi(s).  $\par
 Hence, we need to consider only Gaussian stationary processes. Let us consider the family of examples of real
parts of centered lacunar Gaussian stationary processes:

$$
\xi_{\theta}(t) = \sum_{k=1}^{\infty} k^{-\theta} \ \epsilon_k \ \cos( n(k) t  ), \eqno(4.9)
$$
where $ \{  \epsilon_k, \} \ \{ n(k) \} $ are as in the last subsection,  $ \theta = \const > 1/2. $
 The corresponding  covariation function

$$
R_{\theta}(t,s) = \sum_{k=1}^{\infty} k^{-2 \theta} \ \cos(n(k)t) \ \cos(n(k)s). \eqno(4.10)
$$
 is continuous by virtue  of uniform convergence of the series (4.10).  Let us estimate more precisely
its modulus of continuity. We can use the following known estimate (\cite{Achieser1}, chapter 1:)

$$
|c(k)| \le 0.5 \ \omega(f, \pi/|k|), \  k \ne 0.
$$

 We obtain therefore the following {\it lower} estimation for the modulus of continuity  of the covariation
function $  R: $

 $$
 \omega(R_{\theta}, \pi/n(k)) \ge C \  k^{-2 \theta}, \ k = 1,2,\ldots.
 $$
  Following, the covariation function $ R $ does not satisfy the condition (4.5) for any value $  \theta. $ \par

  But we conclude taking in attention the expression (4.0) based on the properties of lacunar series that the
 r.p. $ \xi_{\theta}(\cdot)  $ is continuous $ (\mod {\bf P} ) $  only iff $ \theta > 1. $\par

\vspace{3mm}

{\bf Remark 4.1.}  Our example (4.9) - (4.10) justifies one of the results of an articles M.Nisio
\cite{Nisio1} and H.Watanabe \cite{Watanabe1}.  Indeed, let $ \eta(t), \ t \in R^d $ be centered separable
Gaussian random field  with covariation function

$$
\rho(t) = {\bf E} \eta(t+s)\eta(s) = \int_{R^d} \cos(\lambda \cdot t) \ F(d \lambda), \eqno(4.11)
$$
where $ F = F(A), \ \subset R^d $ is spectral measure. Denote

$$
B(n) = \{ \lambda, \ |\lambda| \le 2^n \}, \  s_n = F(B(2^{n+1}) - F(B(2^n)).
$$
  Theorem 2 in the article \cite{Watanabe1} asserts that if there exists a monotonically decreasing sequence
$ \{ M(n) \}   $   for which

$$
s(n) \le M(n), \  \sum_{n=1}^{\infty} M^{1/2}(n) < \infty, \eqno(4.12)
$$
then $ \eta(t) $ has continuous sample parts.  Note that for the real part of stationary r.p. $ \xi_{\theta}(t) $ with
$ \theta \in (1/2, 1) $ the condition (4.12) is not satisfied and it is discontinuous and moreover is
unbounded on every non-empty interval.\par

\vspace{4mm}

\section{ Necessary and sufficient conditions for continuity of non-Gaussian processes.}

\vspace{4mm}

 The case of non-Gaussian processes $  \xi(t) $  is more complicated.  Firstly,  the covariation function
for $ \xi(t) $  may do not exists. Secondly, the coefficients in the Karunen-Loev expansion are in general case
dependent. Thirdly, this expression may not convergent uniformly still in the case of continuous process. \par
 So, we consider in this section a separable periodical r.p. $ \xi(t), \ t \in T. $ We can suppose without loss
of generality that it is continuous in probability; as a consequence -  $ \xi(t) $ is measurable over variable $  t $
with probability one. \par

 We need to use the so-called Franklin's  system of continuous functions $  \{ f_k(t) \}, \ k=1,2,\ldots; $
  \cite{Franklin1}, \cite{Ciesielski1}, \cite{Kaczmarz1}. Recall that this system gives an unconditional basis in the
space of all continuous functions. \par

 Denote for the r.p. $ \xi(t) $
$$
\zeta_k = (2 \pi)^{-1} \ \int_T \xi(t) \ f_k(t) \ dt,
$$
so that we can write the formal Fourier - Franklin (FF) expression

$$
\xi(t) = \sum_{k=1}^{\infty} \zeta_k \ f_k(t). \eqno(5.1)
$$

 Let us introduce the following important functionals

 $$
 \beta_n^m = \beta_n^m (\xi) \stackrel{def}{=} \ {\bf  E} \left\{ \arctan || \ \sum_{k=n+1}^m \zeta_k \ f_k \ ||_{\infty} \right\},
 m \ge n+1,  \eqno(5.2.)
 $$
where in the case when $ \beta_n^m $ is not correctly defined, for instance if some coefficient $ \zeta_k $ does
not exists, we put by definition $ \beta_n^m (\xi) = \infty. $  \par

\vspace{4mm}

{\bf Theorem 5.1.} {\it In order to the r.p. } $  \xi(t) $  {\it has continuous with probability one sample path,
is necessary and sufficient that it was stochastic continuous  and that}

$$
\lim_{n \to \infty} \sup_{m \ge n + 1} \beta_n^m (\xi) = 0. \eqno(5.3)
$$

\vspace{4mm}

 {\bf Proof} is alike as in the theorem 3.1. It is sufficient to recall that the functional

 $$
\rho(\xi,\eta) = \rho(\xi-\eta,0) = \rho(\xi -\eta)  =  {\bf E} \ \arctan |\xi - \eta| \eqno(5.4)
 $$
defined on the set of all pair of random variables $ (\xi,\eta) $ is translation invariant distance,
in particular, satisfies the triangle inequality. The set of all r.v. defined on the our probabilistic
space $ \Omega $  equipped  with the distance $ \rho(\cdot, \cdot) $ is complete. \par
 Moreover, the convergence $ \rho(\eta_n, \eta) \to 0, \ n \to \infty $ is equivalent to convergence
 $ \eta_n \to \eta $ in probability. \par

\vspace{3mm}

{\bf Sufficiency.}  Denote the partial sum of FF expression for $ \xi(t) $ as $  \xi_n = \xi_n(t): $

$$
\xi_n(t) = \sum_{k=1}^n  \ \zeta_k \ f_k(t).
$$
 Let $  n(l), \ l=1,2,\ldots $ be strictly increasing  subsequence of natural numbers such that

$$
\beta_{n(l)+1}^{n(l+1) } \le 2^{-l},  \ l \ge l_0 = \const,
$$
and define

$$
S_l = S_l(t) = \sum_{k=n(l) + 1}^{n(l+1)} \ \zeta_k \  f_k(t) = \xi_{n(l+1)}(t) - \xi_{n(l)}(t),
$$
then

$$
{\bf E} \arctan ||S_l||_{\infty} = \beta_{n(l)+1}^{n(l+1) } \le 2^{-l},
$$

$$
 \sum_{l \ge l_0}  {\bf E} \arctan ||S_l||_{\infty} = \sum_{l \ge l_0} \beta_{n(l)+1}^{n(l+1) }  < \infty.
$$
 Therefore, the series

$$
\sum_l \arctan ||S_l||_{\infty}
$$
converges with probability one, and this conclusion also  is true for the series

$$
\sum_l  ||S_l||_{\infty},
$$
following the subsequence of r.p. $ \xi_{n(l)}(t) $  converges uniformly with probability one.
Thus, the r.p. $ \xi(t) $ is continuous $ ( \mod{\bf P}) $. \par

\vspace{3mm}
{\bf  Necessity.} It can be assumed here that the  r.p. $ \xi(t) $ is continuous a.e.   As long as the
Franklin's system $ \{  f_k(\cdot) \} $  formed an unconditional basis in the space of continuous functions
$ C(T), $ the series $ \sum_k \zeta_k f_k(t) $ converges uniformly to the function $ \xi(t)
 ( \mod P): $

$$
{\bf P} ( \lim_{n \to \infty} || \ \xi(\cdot) - \xi_n(\cdot) \ ||_{\infty} = 0) = 1,
$$
all the more

$$
\lim_{n \to \infty} {\bf E} \arctan || \  \xi(\cdot) - \xi_n(\cdot) \ || = 0.
$$

 We get  using the triangle inequality for the $ \rho - $ distance

 $$
 \beta_n^m(R) = {\bf E} \arctan || \xi_m(\cdot) - \xi_n(\cdot) ||_{\infty} \le {\bf E} \arctan || \xi_m(\cdot) - \xi(\cdot) ||_{\infty} +
 $$

$$
{\bf E} \arctan|| \xi(\cdot) - \xi_n(\cdot) ||_{\infty} \to 0, \ n,m \to \infty.
$$
 This completes the proof of theorem 5.1. \par

\vspace{4mm}

\section{ Some sufficient conditions for continuity of non-Gaussian processes.}

\vspace{4mm}

 The author be aware that the conditions of theorem 5.1 are also hard to verify. Further in this section
we will obtain some simple {\it sufficient}  conditions for continuity of  separable stochastic continue periodical non-Gaussian
centered r.p. $ \xi(t)  $ based on some estimations for important for us semi-norm $ \beta_n^m(\xi(\cdot)). $\par
 Our approach here is  development of one  belonging to I.A.Ibragimov \cite{Ibragimov1},  where was applied  the
 approach based on the embedding theorem, which is yet closely related with the approximation theory.\par
We will name the method used in \cite{Ibragimov1} "power method", since it relies on the $ L_p  $ norms of
r.p. $  \xi(t), \xi(t) - \xi(s) $ etc. \par
 We suppose in this section, in contradiction, the existence of {\it exponential} moments of $ \xi(\cdot). $
More exactly, we assume the existence of a so-called {\it generating functional: }

$$
\Phi(\psi) = \Phi_{\xi}(\psi) \stackrel{def}{=} {\bf E} e^{ (\xi,\psi) } = {\bf E} e^{ \int_T \xi(s) \ d \psi(s) }, \eqno(6.1)
$$
where $ \psi(\cdot) $ is deterministic function of bounded variation, $ \psi(0) = 0. $ \par

 Other notations. Let  as before $ n(k), \ k=1,2,\ldots  $ be any  strictly increasing  sequence of natural numbers,

 $$
 Z_k (t)= [V_{n(k+1), p(n(k+1))}* \xi](t) -  [V_{n(k), p(n(k))}* \xi](t) =
 $$
 $$
  [W(n(k+1), n(k),\cdot)*\xi(\cdot)](t) =: [W(\cdot)* \xi(\cdot)](t) ,  \eqno(6.2)
 $$

 $$
 \Psi(\lambda, n(k), n(k+1)) =  (2 \pi)^{-1} \int_T \  {\bf E} e^{ \lambda Z_k(t) } \ dt, \ \lambda = \const > 0; \eqno(6.3)
 $$
 evidently, $ \Psi(\cdot, \cdot, \cdot) $  may be easily expressed through  generating functional $  \Phi(\cdot). $
  Namely, let

 $$
 \psi_t(s) = \int_0^s W(t-x) \ dx,
 $$
 then

 $$
 \Psi(\lambda, n(k), n(k+1)) = (2 \pi)^{-1} \int_T {\bf E} e^{ \lambda \int_T \xi(s) W(t-s) ds } \ dt =
 $$

 $$
 (2 \pi)^{-1} \int_T dt \ {\bf E} \ e^{ \lambda  \ \int_T \xi(s) \ d \psi_t(s)} = (2 \pi)^{-1} \int_T  \ \Phi(\lambda \psi_t(\cdot) ) \ dt.
 $$

 \vspace{4mm}

  It is natural to expect that

  $$
  \lim_{k \to \infty} Z_k (t)= 0
  $$
 and hence

 $$
 \lim_{k \to \infty} \log \Psi(\lambda, n(k), n(k+1)) = 0.
 $$

  Further, define

  $$
  U(n(k), n(k+1)  ) = \inf_{\lambda > 0} \left[  \frac{\log n(k+1)  +  \log \Psi(\lambda, n(k), n(k+1))}{\lambda} \right].\eqno(6.4)
  $$

\vspace{4mm}

{\bf Theorem 6.1.} {\it If for some sequence }  $  \{n(k) \} $

$$
\sum_{k=1}^{\infty} U( n(k), n(k+1) ) < \infty, \eqno(6.5)
$$
{\it then  almost all trajectories of r.p.}  $ \xi(t) $ {\it are continuous.} \par

\vspace{3mm}
{\bf Proof.}\\

{\bf 1. \ Lemma  6.1.}  ( \cite{Ostrovsky209}). Let $  B_n(t), \ t \in T $ be a trigonometrical polynomial of
 degree less or equal than $ n. $ Then

$$
\mes \{t: B_n(t) \ge 0.5 \ || \ B_n \ ||_{\infty} \}  \ge 1/(2n). \eqno(6.6)
$$
\vspace{3mm}
{\bf 2.}  Let us consider the following variable:

$$
I :=  (2 \pi)^{-1} {\bf E} \int_T e^{ \lambda  Z_k(t) } \ dt.
$$
 Since $ Z_k(t) $ is the trigonometrical polynomial  of degree $  \le n(k+1),  $  we can apply  the lemma 6.1:

 $$
 I  \ge (2 \pi)^{-1} {\bf E} \int_{ t: Z_k(t) \ge 0.5 || \ Z_k \ ||_{\infty}  } e^{ \lambda  Z_k(t) } \ dt \ge
 $$

$$
(2 \pi)^{-1} {\bf E} \int_{ t: Z_k(t) \ge 0.5 || \ Z_k \ ||_{\infty}  } e^{ 0.5 \ \lambda  \ || Z_k ||_{\infty} } \ dt \ge
C \ [n(k+1)]^{-1} \  {\bf E} \ e^{0.5 \ \lambda \  || Z_k ||_{\infty} }.
$$
 We apply the Iensen's inequality:

 $$
 {\bf E}  || Z_k ||_{\infty}  \le C  \frac{\log n(k+1) + \log I}{\lambda} = C  \frac{\log n(k+1) + \log \Psi(\lambda, n(k), n(k+1))}{\lambda},
 $$
therefore  $ {\bf E} || Z_k ||_{\infty}  \le $

$$
  C \ \inf_{\lambda > 0}
\left[  \frac{\log n(k+1) + \log \Psi(\lambda, n(k), n(k+1))}{\lambda}  \right] = C \ U(n(k), n(k+1)  ). \eqno(6.7)
$$

\vspace{3mm}
{\bf 3.} We conclude by virtue of condition (6.5) that  the following series  converge:

$$
\sum_{k=1}^{\infty} {\bf E}  || Z_k ||_{\infty}  < \infty,
$$
with him

$$
\sum_{k=1}^{\infty}  || Z_k ||_{\infty}  < \infty \ (\mod {\bf P}).
$$
This completes as before the proof of theorem 6.1. \par

\vspace{4mm}

\section{ Convergence of a  sequence of a random  variables with probability one.}

 \vspace{3mm}
 It is deserve of our attention the investigate  the case of the separable Banach
 space $ c_0 $ of the numerical sequences
$ \{  x(n)  \}, \ n=1,2,\ldots  $  tending to zero  equipped  with the norm

$$
|| \ \{  x(n)  \} \ ||c_0 = \sup_n |x(n)|. \eqno(7.0)
$$
 In detail: let $ \xi =  \{\xi(n) \} $  be a random sequence; find the conditions (necessary conditions and sufficient
conditions) under which  $ {\bf P} (\lim_{n \to \infty} \xi(n) = 0 ) = 1  $ or equally
$ {\bf P}  ( \{ \xi  \} \in c_0 ) = 1.  $\par
  The directly application of theorem  5.1 is very hard; we intend to obtain a simple formulated criterion. \par
We need to introduce a new notations.

$$
\tilde{\xi}(n) :=  \arctan{\xi(n)}, \ \kappa_n^m = \kappa_n^m(\xi) \stackrel{def}{=}
{\bf E} \arctan ( \max_{k=n}^m |\xi(k)|) =
$$

$$
{\bf E}  \max_{k=n}^m |\tilde{\xi}(k)|, \ m \ge n + 1;
$$

$$
A = \{ \omega: \lim_{n \to \infty} \xi(n) = 0   \} = \{ \omega: \lim_{n \to \infty} \tilde{\xi(n)} = 0   \}.
$$

  Further, we set

 $$
 A_Q = \{ \omega: \forall s = 1,2,\ldots \ \exists N = 1,2,\ldots:  \ \max_{ n \in [N, N + Q]} |\tilde{\xi}(n)| < 1/s \},
 $$

 Obviously,

 $$
 {\bf P}(A) = \lim_{Q \to \infty} {\bf P}(A_Q).
 $$

\vspace{4mm}

{\bf Theorem 7.1.} {\it In order to  }  $ {\bf P}(\lim_{n \to \infty} \xi(n) = 0) =1, $ {\it is necessary and sufficient that}

$$
\lim_{n \to \infty} \sup_{m > n} \kappa_n^m (\xi) = 0. \eqno(7.1)
$$

\vspace{4mm}

{\bf Proof. Necessity.}\\

\vspace{3mm}

 Let $ {\bf P}(\lim_{n \to \infty} |\xi(n)| = 0) =1,  $  then with at the same probability

 $$
 \lim_{n \to \infty} \sup_{m > n} | \xi(m)| = 0
 $$
 and  a fortiori

 $$
 \lim_{n \to \infty} \sup_{m > n} \arctan |\xi(m)| = 0
 $$
almost surely. We conclude on the ground  of dominated convergence theorem

$$
\lim_{n \to \infty} {\bf E}  \sup_{m > n} \arctan |\xi(m)| = 0,
$$
which is equivalent to the equality (7.1.) \\

\vspace{4mm}

{\bf Proof. Sufficiency.}\\

\vspace{3mm}

{\bf 1.} Note first of all that

$$
A = \{ \omega: \lim_{n \to \infty} \xi(n) = 0   \} = \{ \omega: \lim_{n \to \infty}  |\tilde{\xi}(n)| = 0   \} =
$$

$$
\cap_s \cup_N  \{  \omega: \sup_{n \ge N} | \tilde{\xi}(n)  |  < 1/s \}  \eqno(7.2)
$$
 and correspondingly

 $$
 A_Q = \cap_s \cup_N  \{  \omega: \max_{n \in [ N, N + Q]} | \tilde{\xi}(n)  |  < 1/s \}.  \eqno(7.3)
 $$

\vspace{3mm}

{\bf 2.} Let the condition (7.1) be satisfied.  We consider a supplementary  events:

$$
B = \overline{A} = \Omega \setminus A, \ B_Q = \overline{A}_Q = \Omega \setminus A_Q.
$$

Elucidation: the set of elementary  events $ B $ may contains (theoretically) also the elementary events when
the limit  does not exists. \par

We can write

 $$
 B =  \cup_s \cap_N \{\omega: \sup_{ n \ge N} | \tilde{\xi}(n)| \ge 1/s  \},
 $$

 $$
 B_Q =  \cup_s \cap_N \{\omega: \max_{ n \in [N, N+Q]} | \tilde{\xi}(n)| \ge 1/s  \} =
 \cup_s C_{s,Q},  \eqno(7.4)
 $$
where

$$
C_{s,Q} = \cap_N \{\omega: \max_{ n \in [N, N+Q]} | \tilde{\xi}(n)| \ge 1/s  \} =  \cap_N D^{(N)}_{s,Q},
$$

$$
D^{(N)}_{s,Q} = \{\omega: \max_{ n \in [N, N+Q]} | \tilde{\xi}(n)| \ge 1/s  \}. \eqno(7.5)
$$

\vspace{3mm}

{\bf 3.} We obtain using the Tchebychev's  inequality:

$$
{\bf P} \left(  D^{(N)}_{s,Q}  \right) \le   \frac{\kappa_N^{N + Q}}{\arctan(1/s)} \to 0, \ N \to \infty,
$$
therefore for all the (natural) values $ s, Q $

$$
 {\bf P} \left( C_{s,Q} \right) =  0,  \eqno(7.6)
$$
following

$$
\forall Q = 1,2,\ldots \ \Rightarrow {\bf P} \left(  B_Q  \right) = 0.
$$

\vspace{3mm}

{\bf 4.} We find

$$
{\bf P} (B) = \lim_{Q \to \infty} {\bf P} \left(B_Q \right)  = 0,
$$
and ultimately $ {\bf P} (A) = 1, $ Q.E.D. \par

 \vspace{3mm}

 {\bf Example 7.1. } Degenerate random processes and sequences.  \par
 \vspace{3mm}
 {\bf Definition 7.1. }  The r.p. $ \xi(t) $ is said to be (linear) degenerate, if

 $$
 \xi(t) = \sum_{k=1}^N \eta_k \ \cdot h_k(t), \ t \in T,  \eqno(7.7)
 $$
 where $ \{ \eta_k \}  $ are linear independent r.v.,  $  \{ h_k(\cdot) \} $  are linear
independent deterministic functions. \par
 We will say in the case $ T = \{  1,2, \ldots \},  $  i.e. when

 $$
 \xi(n) = \sum_{k=1}^N \eta_k \ \cdot h_k(n),   \eqno(7.8)
 $$
 the r.p.  $ \xi(n) $ is called random degenerate sequence. \par
  Let  $ \xi(n) $ be the  random degenerate sequence; then

  $$
  {\bf P} (\xi(n) \to 0) = 1 \ \Leftrightarrow \ \lim_{n \to \infty} {\bf E} \arctan |\xi(n)| = 0. \eqno(7.9)
  $$

\vspace{3mm}

 {\bf Remark 7.1. }

\vspace{3mm}
  We will consider  here the case of the Banach space $ c $ consisting on all the numerical sequences $ \{ x(n) \}  $
with existing the limit

$$
\exists \lim_{n \to \infty} x(n) =: x(\infty)
$$
 with at the same norm as in (7.0). As before, we consider the classical problem:
  let $ \xi =  \{\xi(n) \} $  be a random sequence; find the conditions (necessary conditions and sufficient
conditions) under which  $ {\bf P} ( \exists \lim_{n \to \infty} \xi(n) ) = 1  $ or equally
$ {\bf P}  ( \{ \xi  \} \in c ) = 1.  $\par
 Notations:

$$
\tilde{\xi}(n) :=  \arctan{\xi(n)}, \ \gamma_n^m = \gamma_n^m(\xi) \stackrel{def}{=}
{\bf E} \arctan ( \max_{k=n}^m |\xi(k) - \xi(n)|).
$$

 We find analogously to the theorem 7.1:  $ {\bf P} ( \{\xi \} \in c ) =1 $ if and only if

 $$
 \lim_{n \to \infty} \sup_{m > n} \gamma_n^m (\xi) = 0. \eqno(7.10)
 $$

\vspace{3mm}

 {\bf Remark 7.2. }

\vspace{3mm}

 It is known in the theory of martingales, see e.g. \cite{Burkholder2},  \cite{Burkholder3}, \cite{Hall1},
 chapters 2,3, \cite{Peshkir1} that the estimation of the maximum distribution play a very important role
 for the investigation  of limit theorems, non asymptotical estimations etc. \par
  It follows from our considerations that at the same is true in  more general case of non-martingale
 processes and sequences.

\vspace{3mm}

 {\bf Remark 7.3. } Roughly speaking, the result of theorem 7.1 may be reformulated as follows. Let again
 $ \{ \xi(n) \} $ be a sequence of a r.v., $ \tilde{\xi}(n) = \arctan(|\xi(n)|),  $ and

 $$
 \eta(n) = \sup_{m \ge n} | \tilde{\xi}(n)|. \eqno(7.11)
 $$
  Then

  $$
  \{\omega: \xi(n) \to 0  \}  =  \{\omega: \eta(n) \to 0 \}. \eqno(7.11)
  $$

 But the random sequence $ \{ \eta(n) \} $ is monotonically non-increasing,  therefore the sequence
 $ \{ \eta(n) \}  $ tends to zero {\it with probability one} iff this sequence tends to zero {\it in probability,}
 or equally

 $$
 \lim_{n \to \infty} {\bf E} \eta(n) = 0, \eqno(7.12)
 $$
because the variables $ \eta(n) $ are uniformly bounded. \par

\vspace{4mm}

\section{Concluding remarks.}

 \vspace{4mm}

{\bf 1.  Estimate of modulus of continuity for random process. }\\
\vspace{3mm}

In the  practice, for instance, in the investigation of limits theorem for random processes, may be
appeared the need  of estimation the  modulus of continuity for r. p.  We can use for this purpose
the {\it inverse}  theorems of approximation theory, for example, Stechkin's  estimate (\cite{Timan1},
chapter 6, section 6.1:)

$$
\omega(f, 1/n) \le \frac{c}{n} \cdot  \sum_{m=0}^n E_m(f). \eqno(8.1)
$$
 We ensue for periodical r.p. $ \xi(t): $ \par
\vspace{3mm}
{\bf Proposition 8.1.}

$$
{\bf E} \omega(\xi(\cdot), 1/n) \le \frac{c}{n} \cdot \sum_{m=0}^n {\bf E}  \ E_m(\xi). \eqno(8.2)
$$
\vspace{3mm}

 Note as a slight consequence: \par
\vspace{3mm}
{\bf Proposition 8.2.}\par
\vspace{3mm}
 If $ \lim_{n \to \infty} {\bf E}  \ E_n(\xi) = 0,  $ then

$$
\lim_{n \to \infty} {\bf E} \ \omega(\xi(\cdot), 1/n) = 0
$$
 and therefore

 $$
 \forall \epsilon > 0 \ \Rightarrow  \lim_{\delta \to 0+} {\bf P} ( \omega(\xi(\cdot), \delta) > \epsilon ) = 0.
 $$

\vspace{3mm}

 Moreover,  since the function $ \delta \to \omega(f,\delta) $ is monotonically non-decreasing,

 $$
 {\bf P}( \lim_{\delta \to 0+} \omega(\xi(\cdot), \delta) = 0) = 1. \eqno(8.3)
 $$

\vspace{3mm}

 The last assertion may be strengthened as follows. It is sufficient to suppose instead condition
$ \lim_{n \to \infty} {\bf E}  \ E_n(\xi) = 0  $ to assume
$ \lim_{n \to \infty} {\bf E}  \ \arctan E_n(\xi) = 0;  $ then

$$
\lim_{n \to \infty} {\bf E} \ \arctan \omega(\xi(\cdot), 1/n) = 0;
$$
herewith the conclusion (8.3) remains true. \par

 \vspace{4mm}

{\bf 2.  Application in the theory of limits theorem for random processes. }\\

\vspace{3mm}

 Let $ \xi_{\alpha}(t), \ \alpha \in \cal{A}  $ be a {\it family} of separable periodical processes,
$  \cal{A}  $ be arbitrary set. Assume that for some non-random point $ t_0 \in T $ the family
of one-dimensional r.v. $ \xi_{\alpha}(t_0)  $ is tight.  If

$$
 \lim_{n \to \infty} \sup_{\alpha} {\bf E} \arctan \ E_n(\xi_{\alpha}) = 0, \eqno(8.4)
$$
then the family of distributions $ \mu_{\alpha}(\cdot) $ in the space $ C(T) $ generated by $ \xi_{\alpha}(\cdot):  $

$$
 \mu_{\alpha}(A) = {\bf P} (\xi_{\alpha}(\cdot) \in A )
$$
is weakly  compact. \par

 Indeed, we conclude on the basis of the last subsection

 $$
\lim_{n \to \infty} \sup_{\alpha } {\bf E} \  \arctan  [ \omega(\xi_{\alpha}(\cdot), 1/n) ] = 0.
 $$
 Our proposition follows from theorem 1 in the book  of I.I.Gikhman and A.V.Skorohod
 \cite{Gikhman1},  chapter 9, section 2  after applying the Tchebychev's inequality.\par

 In detail, let us introduce the following  {\it family } of functionals

 $$
 \beta_n^m(\alpha) = \beta_n^m (\xi_{\alpha}), \ \alpha \in \cal{A}, \eqno(8.5)
 $$
see (5.2), (5.3). \par

\vspace{4mm}

{\bf Theorem 8.1.} {\it In order to all the r.p. } $  \xi_{\alpha}(t) $  {\it have continuous with probability one sample path,
and moreover its distributions are weakly compact in the space of continuous  function }  $ C(T), $
{\it is necessary and sufficient that they was stochastic continuous  and that}

$$
\lim_{n \to \infty} \sup_{\alpha} \ \sup_{m \ge n + 1} \beta_n^m (\alpha) = 0. \eqno(8.6)
$$

\vspace{4mm}

Note that the criterion (8.6) is more convenient for applications as in \cite{Gikhman1},  chapter 9, section 2.
 For example, we can use the sufficient condition given by theorem 6.1. \par

\vspace{4mm}

 For instance, let $ \{\alpha \} = 1,2,3,\ldots $ and

 $$
 \zeta_n(t) = n^{-1/2} \sum_{j=1}^n \eta_j(t),
 $$
where $  \{  \eta_j(t) \} $ are continuous independent identical distributed centered r.p. with finite continuous
covariation function  $  R(t,s). $  The weak as  $ n \to \infty $ convergence of the sequence $  \zeta_n(\cdot) $
to the centered continuous Gaussian distributed r.p. $ \zeta_{\infty} $ with at the same covariation function $  R(t,s)  $
in the space $ C(T) $ implies  the Central Limit Theorem (CLT) in this space; see
\cite{Dudley1}, \cite{Kozachenko1}, \cite{Ledoux1}, \cite{Ostrovsky1}; about applications CLT in the space $  C(T) $  in the
Monte-Carlo method see, e.g. \cite{Frolov1},  \cite{Grigorjeva1}.\par

 Note in addition that   the generating functional  for the r.p. $ \zeta_n(\cdot), $ i.e. $ \Phi_{\zeta_n}(\psi) $
may be calculated and uniformly estimated as follows:

$$
  \Phi_{\zeta_n}(\psi) = n \ \Phi_{\eta}(\psi/\sqrt{n}),  \eqno(8.7)
$$

$$
  \Phi_{\zeta_n}(\psi) \le \sup_n \left[ n \ \Phi_{\eta}(\psi/\sqrt{n}) \right] < \infty.  \eqno(8.8)
$$

 \vspace{4mm}

{\bf 3. Many variables. }\\

\vspace{3mm}
 It is no hard to generalize the notions and properties of the degenerate approximation into the functions
of many variables. Besides, we can use instead trigonometrical approximation  the approximation  by means of
algebraic polynomials.\par

\end{document}